\documentclass[reqno,11pt]{amsart}

\usepackage{amssymb,amsfonts, amsmath}
\usepackage{mathrsfs}
\usepackage{color}
\usepackage[normalem]{ulem} 

\newtheorem{thm}{Theorem}[section]
\newtheorem{cor}[thm]{Corollary}
\newtheorem{lem}[thm]{Lemma}
\newtheorem{prop}[thm]{Proposition}
\theoremstyle{definition}

\theoremstyle{remark}
\newtheorem{rem}[thm]{Remark}
\numberwithin{equation}{section}
\newcommand{\Real}{\mathbb R}

\newcommand{\eps}{\varepsilon}
\newcommand{\diag}{\mathrm{diag}}
\newcommand{\F}{\mathscr{F}}

\newcommand{\s}{\mathbb{S}}
\newcommand{\M}{\mathcal{M}}
\newcommand{\one}[1]{\mathbf{1}_{\{#1\}}}

\renewcommand{\P}{\mathsf{P}}
\newcommand{\Q}{\mathrm{Q}}
\newcommand{\E}{\mathsf{E}}
\newcommand{\simplex}{\mathcal{S}^{d-1}}
\DeclareMathOperator*{\argmin}{\mathrm{argmin}}
\DeclareMathOperator*{\argmax}{\mathrm{argmax}}

\newcommand{\lyap}{\lambda}
\newcommand{\Dom}{\mathcal{D}}

\definecolor{gray}{rgb}{0.9,0.9,0.9}

\input cyracc.def

\begin{document}

\title{An ergodic theorem for filtering with applications to stability}%
\author{Pavel Chigansky}
\address{Department of Mathematics, The Weizmann Institute of Science, Rehovot 76100, Israel}
\email{pavel.chigansky@weizmann.ac.il}

\thanks{Research supported by a grant from the Israel Science Foundation}

\date{6, July, 2006}%
\begin{abstract}
Ergodic properties of the signal-filtering pair are studied for
continuous time finite Markov chains, observed in white noise. The
obtained law of large numbers is applied to the stability problem
of the nonlinear filter with respect to initial conditions. The
Furstenberg-Khasminskii formula is derived for the top Lyapunov
exponent of the Zakai equation and is used to estimate the
stability index of the filter.
\end{abstract}
\maketitle

\section{Introduction}

Consider a pair of continuous time random processes $(X,Y)=(X_t,Y_t)_{t\ge 0}$, where the {\em signal}
component $X$ is a Markov chain, taking values in a finite alphabet $\s=\{a_1,...,a_d\}$, with transition intensities
matrix $\Lambda$ and initial distribution $\nu$. The {\em observation} process $Y$ is given by
\begin{equation}
\label{Y}
Y_t = \int_0^t h(X_s)ds +\sigma B_t,\quad t\ge 0
\end{equation}
with an $\s\mapsto \Real$ function $h$, constant $\sigma>0$ and a Brownian motion $B=(B_t)_{t\ge 0}$, independent of $X$.

The filtering problem is to calculate the conditional probabilities $\pi_t(i)=\P(X_t=a_i|\F^Y_t)$, $i=1,...,d$
where $\F^Y_t=\sigma\{Y_s, s\le t\}$, which are the main building blocks of the optimal MSE and MAP signal estimates
$$
\hat X^{\mathrm{mse}}_t = \sum_{i=1}^d a_i\pi_t(i)
\quad \text{and}\quad
\hat X^{\mathrm{map}}_t = \argmax_{a_i\in \s} \pi_t(i)
$$
given the trajectory of $Y$ up to time $t$.

The vector $\pi_t$ satisfies the It\^o stochastic differential equation (SDE) (\cite{W}, see also \cite{LS})
\begin{equation}
\label{W}
d\pi_t = \Lambda^* \pi_tdt + \sigma^{-2}\big(\diag(\pi_t)-\pi_t\pi^*_t\big)h\big(dY_t-\pi_t(h) dt\big), \quad \pi_0=\nu,
\end{equation}
where $h$ stands for the column vector with entries $h_i:=h(a_i)$, $i=1,...,d$. Hereafter the following notations are
used: $\diag(x)$, $x\in\Real^d$ stands for a scalar matrix with entries $x_i$
and $x^*$ is transposed of $x$. The space of probability measures on $\s$ is identified with the simplex
$\simplex=\{x\in\Real^d: x_i\ge 0, \sum_{i=1}^d x_i =1\}$ and $
\mu(f):=\sum_{i=1}^d \mu_if(a_i)
$
for $f:\s\to\Real$ and $\mu\in\simplex$.
For vectors and matrices $|\cdot|$ denotes the $\ell_1$-norm, i.e. $|x|=\sum_{ij} |x_{ij}|$.
Finally $\F^X_{s,t}:=\sigma\{X_r, s\le r\le t\}$, $\F^Y_{s,t}:=\sigma\{Y_r-Y_s, s\le r\le t\}$ and
$\F_{s,t}:=\F^Y_{s,t}\vee \F^X_{s,t}$ ($\F^X_t:=\F^X_{0,t}$, etc. are written for brevity).
As usually all the statements involving random objects are understood to hold $\P$-a.s.

The process $\bar B_t := \sigma^{-1}\int_0^t \big(dY_s - \pi_s(h)ds\big)$ is the {\em innovation} Brownian
motion with respect to the filtration $\F^Y_t$ and so the solution of \eqref{W} is a Markov process. In this
paper we deal with the ergodic properties of the process $\pi=(\pi_t)_{t\ge 0}$ (and more generally of the pair $(X,\pi)$).
It should be noted that the Wonham SDE \eqref{W} does not fit the standard ergodic theory of diffusions (see e.g. \cite{Kh}),
which usually requires either uniform non-degeneracy of the diffusion matrix or Hormander's hypoellipticity conditions.
The study of ergodic properties of  the filtering process for general Markov signals was initiated in early 70's by H.Kunita in \cite{Ku},
but one of the main arguments in his results appears to have a serious gap, leaving the problem open (see \cite{BCL} for further details).

Our main motivation for studying the ergodic properties of the filtering process $\pi$, is the stability of the filtering
equation \eqref{W} with respect to the initial condition. It is not hard to see (e.g. as in \cite{RS}) that \eqref{W} has a unique strong solution, if started
from any $\bar\nu\in\simplex$, possibly different from $\nu$. Denote the corresponding solution by $\bar\pi=(\bar \pi_t)_{t\ge 0}$.
The filter is said to be {\em asymptotically stable}, if
\begin{equation}
\label{stab}
\lim_{t\to\infty}|\pi_t-\bar{\pi}_t|=0.
\end{equation}
The problem of stability is to find the conditions in terms of  $\Lambda$, $h$, $\sigma$ and $(\nu,\bar\nu)$,
which guarantee \eqref{stab}.

There has been much progress in solving the stability problem during the past decade for various filtering models: we mention \cite{AZ97a,AZ97b,BO,LO,DG,BCL}
for a few.
In particular, \eqref{stab} for the Wonham SDE has been verified recently in \cite{BCL}: it turns out that a stronger exponential convergence
\begin{equation}
\label{exstab}
\gamma:=\lim_{t\to\infty}\frac{1}{t}\log|\pi_t-\bar{\pi}_t|<0,
\end{equation}
holds for any $h$, $\sigma>0$ and $(\nu,\bar \nu)$, if the chain $X$ is ergodic, i.e. the limit probabilities $\mu_i:=\lim_{t\to\infty}\P(X_t=a_i)$, $i=1,...,d$
exist, are positive and independent of $\nu$. The {\em stability index} $\gamma$ is of significant practical value, as it quantifies the rate of
convergence in \eqref{stab}. Regretfully it is quite elusive for explicit calculation (see Section \ref{sec-apl} for a summary of the available
bounds for $\gamma$).

The main result of this paper is the following ergodic theorem for the signal-filtering pair $(X,\pi)$
\begin{thm}\label{e-thm}
Assume that $X$ is ergodic, then the Markov-Feller process $(X,\pi)$ has a unique invariant measure $\M$, such that
for any continuous $g$
\begin{multline}
\label{raz}
\lim_{t\to\infty}\frac{1}{t}\int_0^t g(X_s,\pi_s)ds=\int_{\s\times\simplex}g(x,u)\M(dx,du)= \\
\int_{\simplex} \sum_{i=1}^d u_i g(a_i, u)\M_\pi(du)=\lim_{t\to\infty} \E g(X_t,\pi_t),
\end{multline}
where $\M_\pi$ is the $\pi$-marginal of $\M$.
\end{thm}

The ergodic properties \eqref{raz} are derived from the stability \eqref{stab} and turn to be
useful in the study of exponential convergence in \eqref{exstab} within the framework due to H.Furstenberg and R.Khasminskii
(see \cite{Kh}, \cite{Arnold}). Besides providing an additional insight into the problem, the method simplifies derivation of several
already known upper bounds for $\gamma$ (\cite{AZ97a}, \cite{DZ}). Also it allows to obtain a closed form formula
and some higher order asymptotic expansions for $\gamma$ in the two-dimensional  case $d=2$.

The proof of Theorem \ref{e-thm} is given in the next Section \ref{e-pr}. In Section \ref{sec-apl} the Law of Large Numbers (LLN) from
\eqref{raz} is applied to the problem of stability, using the Lyapunov exponents technique introduced in \cite{AZ97a}.

\section{The proof of Theorem \ref{e-thm}}\label{e-pr}

Since $X$ is a Markov process and $\pi$ is the unique strong solution of the time homogeneous diffusion SDE \eqref{W},
driven by $X$ and an independent Brownian motion $B$, the pair $(X,\pi)$ is Markov. The Feller property for $\pi$, verified in Theorem 2.3, \cite{Ku},
is inherited by the pair $(X,\pi)$. Being a Feller-Markov process on a compact state space, the transition probability of $(X,\pi)$
has at least one invariant measure (see e.g. proof of Theorem 3.1, \cite{Ku}).

It will be convenient to consider the {\em semiflow} generated by \eqref{W}, i.e. the family of random maps
$u\mapsto\pi_{s,t}(u)$, $t\ge s\ge 0$ of $\simplex$ into itself, obtained by solving \eqref{W} on $[s,t]$ subject to $\pi_s=u$.
In fact it can be seen, that $\P$-a.s. $u\mapsto\pi_{s,t}(u)$ is a smooth injective $\simplex\to\simplex$ map,
satisfying $\pi_{r,t}(\pi_{s,r}(u))=\pi_{s,t}(u)$, $t\ge r\ge s\ge 0$ (see e.g. \cite{Kuflows}).
For example the process $\bar\pi_t$, defined in the Introduction, is nothing but $\pi_{0,t}(\bar\nu)$
and we will use both notations when no confusion occurs.

The following result is a straightforward modification of Theorem 4.1 in \cite{BCL}, whose proof is
omitted
\begin{thm}
Assume that $X$ is ergodic, then for any $\nu\in\simplex$, $s\ge 0$ and $\F_s$-measurable random vectors $u,v\in\simplex$
\begin{equation}
\label{rands}
\varlimsup_{t\to\infty}\frac{1}{t}\log\big|\pi_{s,t}(u)-\pi_{s,t}(v)\big|<0.
\end{equation}
\end{thm}
The uniqueness of the invariant measure for $(X,\pi)$ is established in Theorem 7.1, \cite{BK} under assumption \eqref{rands}.

We will use the following construction for the Markov chain $X$. Let $\eta=(\eta_t)_{t\ge 0}$ be the solution of the It\^o
equation
\begin{equation}
\label{ito}
\eta_t = \eta_0 + \int_0^t dN^*_s \eta_{s-},
\end{equation}
where $N_s$ is a matrix, whose off-diagonal elements are independent Poisson processes with intensities $\lambda_{ij}$ and
$N_s^{ii}=-\sum_{j\ne i}N_s^{ij}$. If $\eta_0$ is a random vector, taking values in the standard basis $\mathcal{E}:=\{e_1,...,e_d\}$ of $\Real^d$, with
$\P(\eta_0=e_i)=\nu_i$ and independent of $N$, then $\eta_t\in \mathcal{E}$ for all $t\ge 0$ and the random process $X_t=\sum_{i=1}^d a_i \eta_t(i)$ is
a Markov chain\footnotemark with transition rates matrix $\Lambda$ and initial distribution $\nu$.
\footnotetext{without loss of generality $a_i\ne a_j$, $i\ne j$ can be assumed}

Let $(\tilde X_0, \tilde \pi_0)$ be a random variable with distribution $\M$. Let $\tilde X$ be the stationary Markov chain, e.g.
generated by \eqref{ito}, and $\tilde Y$ be defined by \eqref{Y} with $X$ replaced with $\tilde X$. Finally let $\tilde \pi$ be
the solution of \eqref{W}, driven by $\tilde Y$ and started from $\tilde\pi_0$. Evidently the process $(\tilde X, \tilde \pi)$ is
stationary. It is well known that the set of invariant measures of a stationary Feller-Markov process is closed and convex and that
the extremal measures are ergodic (see e.g. Section 6.3 \cite{Var}). This means that if an invariant measure is unique, it is necessarily ergodic
and hence the first equality in \eqref{raz} holds for $(\tilde X,\tilde \pi)$ by the Birkhoff-Kintchine  LLN
for any $\M$-integrable function $g$.
The second equality follows from the property of conditional expectations:
\begin{multline*}
\E g(\tilde X_t,\tilde \pi_t) = \E \sum_{i=1}^d \one{\tilde X_t=a_i}g(a_i,\tilde \pi_t) = \\
 \E \sum_{i=1}^d \E\big(\one{\tilde X_t=a_i}\big|\F^{\tilde Y}_t\vee \tilde \pi_0\big)g(a_i,\tilde \pi_t) =
\E \sum_{i=1}^d \tilde \pi_t(i) g(a_i,\tilde \pi_t).
\end{multline*}

To verify \eqref{raz} for  $(X,\pi)$, we will use the following coupling
\begin{lem}\label{c-lem}
Assume that $\Lambda$ corresponds to an ergodic chain and define
(with the usual convention $\inf\{\emptyset\}=\infty$)
$$
\tau=\inf\{t\ge 0: X_t = \tilde X_t\},
$$
where $X_t$ and $\tilde X_t$ are Markov chains corresponding to the solutions of \eqref{ito}, started from independent $\eta_0$ and $\tilde \eta_0$, independent
of $N$. Then $\lim_{n\to\infty}\one{\tau\ge n}=0$.
\end{lem}
\begin{proof}
Consider the embedded discrete time Markov chain $Z_n := X_n$, $n\in\mathbb{Z}_+$ with the transition probabilities matrix
$G:=\exp(\Lambda)$ and initial distribution $\nu$. Since $X$ is ergodic, all the entries of $G$ are positive (see e.g. \cite{Nor}).
Similarly  $\tilde Z_n:=\tilde X_n$, $n\in\mathbb{Z}_+$ is a Markov chain with the same transition matrix $G$ and initial distribution $\mu$.
Moreover the pair $(Z_n,\tilde Z_n)$ is a Markov chain as well. Hence, on the set $\{Z_{n-1}\ne \tilde Z_{n-1}\}$
$$
\P\big(Z_n\ne \tilde Z_n\big |\F^Z_{n-1}\vee \F^{\tilde Z}_{n-1}\big) \le 1-\sum_{i=1}^d \min_{k\ne \ell}G_{ki}G_{\ell i}=:r<1,
$$
and
\begin{align*}
\P(\tau \ge n) =& \P\big(X_t\ne \tilde X_t, \forall t\le n\big)\le \P\big(Z_k\ne \tilde Z_k, \forall k\le n\big)=\\
& \E \prod_{k=0}^{n-1} \one{Z_k\ne \tilde Z_k} \P\big(Z_n\ne \tilde Z_n\big |\F^Z_{n-1}\vee \F^{\tilde Z}_{n-1}\big) \le \\
& r \E \prod_{k=0}^{n-2} \one{Z_k\ne \tilde Z_k} \P\big(Z_{n-1}\ne \tilde Z_{n-1}\big |\F^Z_{n-2}\vee \F^{\tilde Z}_{n-2}\big) \le ...\le Cr^n,
\end{align*}
with a positive constant $C$. The required claim holds by the Borel-Cantelli Lemma. \qed
\end{proof}

Suppose that $X$ and $\tilde X$ are as defined in Lemma \ref{c-lem} and $\pi$ and $\tilde \pi$ be the corresponding
filtering processes, generated by \eqref{W}, driven by $Y$ and $\tilde Y$ respectively. Note that on the set $\{\tau < n\}$
the increments of $Y$ and $\tilde Y$ coincide after time $n$, i.e. $Y_t-Y_n=\tilde Y_t-\tilde Y_n$, for $t\ge n$ and hence on this set
$$
\tilde\pi_t =  \tilde \pi_{n,t}\big(\tilde\pi_{0,n}(\tilde \pi_0)\big)=\pi_{n,t}\big(\tilde\pi_{0,n}(\tilde \pi_0)\big).
$$
Then by \eqref{rands}
\begin{multline*}
|\pi_t-\tilde \pi_t| = |\pi_{0,t}(\nu)-\tilde \pi_{0,t}(\tilde \pi_0)| = \\
 |\pi_{0,t}(\nu)-\tilde \pi_{0,t}(\tilde \pi_0)|\one{\tau \ge n}
+\big|\pi_{n,t}\big(\pi_{0,n}(\nu)\big)-\tilde \pi_{n,t}\big(\tilde \pi_{0,n}(\tilde\pi_0)\big)\big|\one{\tau <  n} \le\\
  2\one{\tau \ge n}
+\big| \pi_{n,t}\big(\pi_{0,n}(\nu)\big)- \pi_{n,t}\big(\tilde \pi_{0,n}(\tilde\pi_0)\big)\big|\xrightarrow{t\to\infty} 2\one{\tau \ge n}
\end{multline*}
for any fixed $n\ge 0$. By Lemma \ref{c-lem} and arbitrariness of $n$, we have\footnote
{
Note that unlike in \eqref{pi-pi}, the filtering processes in \eqref{stab} are
generated by \eqref{W}, driven by {\em the same} observations.
}
\begin{equation}\label{pi-pi}
\lim_{t\to\infty}|\pi_t-\tilde \pi_t|=0.
\end{equation}
Now the first equality in \eqref{raz} follows from  Lemma \ref{c-lem} and \eqref{pi-pi}, the LLN for $(\tilde X,\tilde \pi)$ and the fact, that $\P$-a.s. convergence implies
$\P$-a.s. Cesaro convergence. The last equality in \eqref{raz} is obtained by means of triangle inequality and dominated convergence. \qed
\section{Application to filter stability}\label{sec-apl}

\subsection{A brief survey}
As was already mentioned in the Introduction, the stability of \eqref{W} for ergodic chains is always exponential in the sense that the limit in
\eqref{exstab} exists and is strictly negative. The existence of the limit $\gamma$ follows from the Oseledec Multiplicative Ergodic Theorem (MET) and, moreover,
it may take no more than a finite number of values, depending\footnote{as will be clarified below, $\gamma$ turns to be independent of $(\nu,\bar{\nu})$
for $d=2$; the actual dependence of $\gamma$  on $(\nu,\bar\nu)$ in the case $d>2$ remains unclear.}
on $(\nu,\bar\nu)$ (see \cite{AZ97a} for details). Unfortunately it is extremely hard to come up with the exact expression for $\gamma$ and usually one
gets only qualitative information on its dependence on the model parameters.

The first result of this type was obtained in \cite{DZ}, where the following local {\em low signal-to-noise} asymptotic has been derived
\begin{equation}
\label{lyap}
\varlimsup_{\sigma\to \infty}\varlimsup_{t\to\infty}\frac{1}{t}\log \Big(\lim_{\bar{\nu}\to \nu}\frac{|\pi_t-\bar{\pi}_t|}{|\nu-\bar{\nu}|}\Big) \le
\gamma_{\max}(\Lambda),
\end{equation}
where $\gamma_{\max}(\Lambda)$ is the spectral gap of $\Lambda$, i.e. the largest non-zero real part of the eigenvalues of $\Lambda$. Roughly speaking \eqref{lyap}
means that in the low signal-to-noise regime the filter is at least as stable as the chain itself, if  $\bar{\nu}$ and $\nu$ are close enough.

The following estimates of $\gamma$ have been derived in \cite{AZ97a}
\begin{align}
& \gamma \le -2\min_{i\ne j}\sqrt{\lambda_{ij}\lambda_{ji}} \label{AZ}\\
& \varlimsup_{\sigma\to 0} \sigma^2 \gamma(\sigma) \le -\frac{1}{2} \sum_{j=1}^d \mu_j\min_{j\ne i}\big(h_j-h_i\big)^2\label{AZu}\\
& \varliminf_{\sigma\to 0} \sigma^2 \gamma(\sigma) \ge - \frac{1}{2} \sum_{j=1}^d \mu_j \sum_{i=1}^d \big(h_j-h_i\big)^2\label{AZl},
\end{align}
where $\gamma(\sigma)$ is written  to emphasize the dependence on the noise intensity $\sigma$ and $\mu$ is the unique
invariant measure of $X$. Note that the bound \eqref{AZ} is independent of the observation
parameters $h$ and $\sigma$, which suggests that to a certain extent stability is inherited by the filter from the signal itself.
Namely the right hand side of \eqref{AZ} remains  negative if all the transition intensities of $X$ are strictly positive. Clearly this
condition is much stronger than ergodicity.
In contrast to \eqref{AZ} the bounds \eqref{AZu} and \eqref{AZl} reveal  the dependence of $\gamma$ on the noise intensity in the
{\em high signal-to-noise} regime: the stability index
decreases  quadratically with $\sigma$ if the image of $\s$ under $h$ has at least one unique point. This perfectly agrees with the fact that
otherwise $\gamma(\sigma)$ may converge to zero as $\sigma\to 0$ as mentioned in \cite{DZ}.

Recently a non asymptotic version of \eqref{AZ} was derived in \cite{BCL}
\begin{equation}
\label{BCL}
|\pi_t-\bar{\pi}_t|\le C\exp\Big(-2t \min_{i\ne j}\sqrt{\lambda_{ij}\lambda_{ji}}\Big),
\end{equation}
where $C>0$ depends only on $(\nu,\bar{\nu})$ and, moreover,
$$
\gamma \le - \sum_{i=1}^d \mu_i \min_{j\ne i}\lambda_{ij}.
$$
Unlike \eqref{AZ} or \eqref{BCL}, this bound remains nontrivial as long as the signal is ergodic and $\Lambda$ has at least
one row with all non-zero entries.

\subsection{The method of Lyapunov exponents \cite{AZ97a}}

In this section we sketch the main idea of the approach to filter stability due to  R.Atar and O.Zeitouni \cite{AZ97a}.
Recall that the solutions of SDE \eqref{W} $\pi_t$ and $\bar{\pi}_t$ can be obtained by solving the
linear Zakai SDE
\begin{equation}
\label{zak}
d\rho_t = \Lambda^* \rho_t dt + \sigma^{-2}\diag(h)\rho_t dY_t,
\end{equation}
and  normalizing
$
\pi_t = \rho_t/|\rho_t|
$
and
$
\bar{\pi}_t = \bar{\rho}_t/|\bar{\rho}_t|,
$
where $\rho_t$ and $\bar{\rho}_t$ denote the solutions of \eqref{zak} subject to $\rho_0=\nu$ and $\rho_0=\bar{\nu}$ respectively.

For a pair of vectors $x,y\in\Real^d$, let $x\wedge y$ denote the exterior product, which can be identified with the
matrix $(x_iy_j-x_jy_i)$, $i\le i,j\le d$. Then the following elementary inequalities hold
$$
\frac{1}{2}\frac{|\rho_t\wedge \bar{\rho}_t|}{|\rho_t||\bar{\rho}_t|}
\le |\pi_t-\bar{\pi}_t|\le
\frac{|\rho_t\wedge \bar{\rho}_t|}{|\rho_t||\bar{\rho}_t|},
$$
which agrees with the fact that the distance between vectors in $\simplex$ can be measured by the angle they form.
These inequalities suggest that
\begin{equation}
\label{Lyap}
\lim_{t\to\infty}\frac{1}{t} \log |\pi_t-\bar{\pi}_t| =
\lim_{t\to\infty}\frac{1}{t} \log |\rho_t\wedge \bar{\rho}_t|-
\lim_{t\to\infty}\frac{1}{t} \log |\rho_t| -
\lim_{t\to\infty}\frac{1}{t} \log |\bar{\rho}_t|
\end{equation}
if the limits in the right hand side exist. Verifying the assumptions of MET (see e.g. \cite{Arnold}),
it was shown in \cite{AZ97a} that for any $\nu$ and $\bar{\nu}$ the limits exist
\begin{equation}\label{top}
\lyap_1 =: \lim_{t\to\infty}\frac{1}{t} \log |\rho_t| = \lim_{t\to\infty}\frac{1}{t} \log |\bar{\rho}_t|,
\end{equation}
are non-random, do not depend on $(\nu,\bar{\nu})$ and equal to the top Lyapunov exponent of the equation \eqref{zak}.
Notably any solution of \eqref{zak} started from a vector from $\simplex$ ``picks up'' the top Lyapunov exponent.
This is a consequence of the Perron-Frobenious theorem applied in \cite{AZ97a} to the positive stochastic flow generated by \eqref{zak}.
The matrix $\rho_t\wedge \bar{\rho}_t$ satisfies a linear SDE (see \eqref{zeq} below) as well and thus is in the scope of MET:
\begin{equation}\label{pretop}
\lim_{t\to\infty}\frac{1}{t} \log |\rho_t\wedge \bar{\rho}_t| \le \lyap_1+\lyap_2,
\end{equation}
where $\lyap_2$ is the second Lyapunov exponent of \eqref{zak}. Roughly speaking \eqref{pretop} means that the area between the two solutions
of \eqref{zak} grows exponentially with a rate not exceeding the sum of two Lyapunov exponents, just as in the case of deterministic linear equations with
constant coefficients. Unlike in \eqref{top}, only inequality can be claimed in \eqref{pretop}
since the flow $\rho_t\wedge \bar{\rho}_t$ in not positive anymore and different initial conditions may in principle correspond to different
Lyapunov exponents.

Assembling \eqref{Lyap}, \eqref{top} and \eqref{pretop} together it is concluded in \cite{AZ97a} that the stability of \eqref{W} is controlled by the
Lyapunov spectral gap of \eqref{zak}:
\begin{equation}\label{lyaplyap}
\lim_{t\to\infty}\frac{1}{t} \log |\pi_t-\bar{\pi}_t|\le \lyap_2-\lyap_1\le 0.
\end{equation}
Though conceptually appealing, the latter does not immediately provide an easy way to verify the desired stability, since
$\lyap_1$ and $\lyap_1+\lyap_2$ are usually not easy to calculate.

Nevertheless the Lyapunov exponents turn to be amenable to asymptotic expansions in terms of parameter $\sigma$: the bound \eqref{AZu} is a
combination of the estimates
\begin{equation}
\label{est1}
\lim_{\sigma\to 0}\sigma^2\lyap_1(\sigma) = \frac{1}{2} \mu (h^2)
\end{equation}
and
\begin{equation}
\label{est12}
\varlimsup_{\sigma\to 0}\sigma^2\big(\lyap_1(\sigma)+\lyap_2(\sigma)\big) \le \frac{1}{2} \mu(h^2) + \mu(h h_{\mathrm{nbr}})-\frac{1}{2}
\mu(h^2_{\mathrm{nbr}}),
\end{equation}
where  $h_{\mathrm{nbr}}(a_i):=h(\mathrm{nbr}(a_i))$ with $\mathrm{nbr}(a_i)=\argmin_{a_j\ne a_i}|h(a_i)-h(a_j)|$.
Both \eqref{est1} and \eqref{est12} are obtained in \cite{AZ97a} using the Feynman-Kac type formulae adapted to the filtering context.

\subsection{Theorem \ref{e-thm} and Furstenberg-Khasminskii formulae}

The objective of this section is to show how the Lyapunov exponents of \eqref{zak} can be
estimated by means of so called Furstenberg-Khasminskii formulae.
\begin{prop}
Assume that $X$ is ergodic, then
\begin{equation}
\label{L1}
\lyap_1 =\frac{1}{2}\sigma^{-2} \int_{\simplex} \big(u(h)\big)^2 \M_\pi(du),
\end{equation}
where $\M_\pi$ is the invariant measure of $\pi$ from Theorem \ref{e-thm}.
\end{prop}

\begin{proof}

Due to \eqref{top}, we can work with $\rho_t$. Since the entries of $\rho_t$ are nonnegative,  $|\rho_t|=\sum_{i=1}^d \rho_t(i)$
and by the Ito formula
\begin{multline*}
d\log |\rho_t| = \frac{1}{|\rho_t|}\sum_{i=1}^d d\rho_t(i) - \frac{1}{2}\frac{1}{|\rho_t|^2} \sigma^{-2}\Big(\sum_{i=1}^d h_i \rho_t(i)\Big)^2dt\stackrel{\dagger}{=}\\
\sigma^{-2}\sum_{i=1}^d h_i \pi_t(i) dY_t - \frac{1}{2} \sigma^{-2}\Big(\sum_{i=1}^d h_i \pi_t(i)\Big)^2dt,
\end{multline*}
where the property $\sum_{j=1}^d \lambda_{ij}=0$ was used in $\dagger$. The idea of studying growth rate of the solutions of linear
SDEs by projecting them onto the unit sphere ($\simplex$ in this case) and averaging with respect to the invariant measure of the
``{\em angle}'' process $\pi_t$, dates back to the works of  H.Furstenberg and R. Khasminskii (see \cite{Kh}) and today constitutes an important
part of the theory of random dynamical systems (see e.g. \cite{Arnold}).

Since $\pi_t(h)$ is  bounded, $\lim_{t\to\infty}\frac{1}{t}\int_0^t \pi_s(h)dB_s=0$ (see e.g. Lemma 7.1, Chapter VI \cite{Kh}),
and the required \eqref{L1} holds by Theorem \ref{e-thm}:
\begin{multline*}
\lyap_1 = \lim_{t\to\infty}\frac{1}{t} \log |\rho_t| =
\lim_{t\to\infty} \frac{1}{t}\int_0^t \Big(
\sigma^{-2}\pi_s(h) dY_s - \frac{1}{2} \sigma^{-2}\big(\pi_s(h)\big)^2ds\Big)=\\
\lim_{t\to\infty} \frac{1}{t}\int_0^t \sigma^{-2} \Big(
\pi_s(h) h(X_s) - \frac{1}{2} \big(\pi_s(h)\big)^2\Big)ds = \frac{1}{2}\sigma^{-2}\int_{\simplex} \big(u(h)\big)^2 \M_\pi(du).
\end{multline*}
\qed
\end{proof}

\begin{rem}\label{rem-1}
Note that \eqref{est1} easily follows from \eqref{L1}, since (the superscript $\sigma$ is used to emphasize the dependence on $\sigma$)
$$\E\big(\tilde \pi^\sigma_t(h)\big)^2 = \E h^2(\tilde X_t) - \E\big(h(\tilde X_t)-\tilde \pi^\sigma_t(h)\big)^2,$$ and
$
\lim_{\sigma\to 0}\E\big(h(\tilde X_t)-\tilde \pi^\sigma_t(h)\big)^2=0,
$
$t>0$.
In fact \eqref{L1} provides even more information, since $\M^\sigma_\pi$, though not computable explicitly in general, enjoys nice concentration properties
revealed by R.Khasminskii and O.Zeitouni in \cite{KZ} and  G. Golubev in \cite{Gol}. Consider the slow stationary Markov chain $\tilde X^\eps$
with generator $\eps\Lambda$, and let  $\tilde Y^\eps_t$ satisfy \eqref{Y} with $X$ replaced by $\tilde X^\eps$ and $\tilde \pi^\eps_t$ be the vector of
corresponding conditional probabilities, i.e. the solution of \eqref{W} with $\Lambda$ and $Y$ replaced by $\eps\Lambda$ and $\tilde Y^\eps$ respectively.
Then, assuming that all $h_i$'s are different and using Theorem \ref{e-thm} and the continuous time analog of  Theorem 1 in \cite{Gol}, one obtains
the asymptotic $\eps\to 0$
$$
\begin{aligned}
&\E\big(h(\tilde X^\eps_t)-\tilde \pi^\eps_t(h)\big)^2 = \\
&\hskip 0.8in\big(1+o(1)\big)\eps\log\eps^{-1}\sum_{i=1}^d \sum_{j\ne i} \frac{2\mu_i\lambda_{ij}}{(h_i-h_j)^2}(h_j-h_j)^2 = \\
&\hskip 0.8in\big(1+o(1)\big)\eps\log\eps^{-1}\sum_{i=1}^d 2\mu_i \sum_{j\ne i}\lambda_{ij}=
\big(1+o(1)\big)\eps\log\eps^{-1}\sum_{i=1}^d 2\mu_i |\lambda_{ii}|.
\end{aligned}
$$
By a time scaling or directly from the Kolmogorov-Fokker-Plank equation for the density of $\M_\pi$, the latter implies
$$
\E\big(h(\tilde X_t)-\tilde \pi^\sigma_t(h)\big)^2 = \big(1+o(1)\big)\sigma^2\log\sigma^{-2}\sum_{i=1}^d 2\mu_i |\lambda_{ii}|, \quad \sigma\to 0.
$$
and in turn
\begin{equation}
\label{refined}
\lyap_1(\sigma)=\sigma^{-2}\frac{1}{2}\mu(h^2) -
\big(1+o(1)\big)\log\sigma^{-2}\sum_{i=1}^d \mu_i |\lambda_{ii}|, \quad \sigma\to 0.
\end{equation}
Note that the second term in the expansion of $\lyap_1(\sigma)$ as $\sigma\to 0$ does not depend on $h$
and is negative, so that the top Lyapunov exponent is actually ``slightly smaller'' than its limit \eqref{est1}.
\end{rem}

\begin{rem}
If $\sigma\to \infty$, the invariant measure $\M^\sigma_\pi$ concentrates at $\delta_{\mu}$ and hence
by \eqref{L1}
\begin{equation}\label{cr}
\lim_{\sigma\to \infty}\sigma^2\lyap_1(\sigma) = \frac{1}{2}\big(\mu(h)\big)^2.
\end{equation}
Indeed the stationary process $\tilde \pi^\sigma$ satisfies
$$
\tilde\pi^\sigma_t = e^{\Lambda^* t}\tilde \pi^\sigma_0 + \sigma^{-1} \int_0^t e^{\Lambda^*(t-s)} \big(\diag(\tilde\pi^\sigma_t)-\tilde\pi^\sigma_t\tilde\pi^{\sigma*}_t\big)h d\tilde B_t,
$$
where $\tilde B$ is the innovation Brownian motion $\tilde B_t = \sigma^{-1} \int_0^t \big(dY_s-\tilde \pi^\sigma_s(h)ds\big)$. Hence,
for any fixed $t>0$,
$
\lim_{\sigma\to \infty}\tilde \pi^\sigma_t = e^{\Lambda^* t}\tilde \pi_0.
$
But since the chain $X$ is ergodic, $\lim_{t\to\infty}e^{\Lambda^* t}\tilde \pi_0= \mu$, and hence
$\M^\sigma_\pi$ converges weakly to $\mu$ as $\sigma\to \infty$.

In fact \eqref{cr} can be further refined, using the results from \cite{Ch3}:
\begin{equation}
\label{crref}
\lim_{\sigma\to \infty}\lyap_1(\sigma) = \frac{1}{2}\big(\mu(h)\big)^2\sigma^{-2} + \frac{1}{2} h^* \Gamma h\sigma^{-4}  + o(\sigma^{-4}),
\end{equation}
where $\Gamma$ is the unique solution of the algebraic Lyapunov equation
$$
0=\Lambda^* \Gamma + \Gamma\Lambda + \big(\diag(\mu)-\mu\mu^*\big)hh^*\big(\diag(\mu)-\mu\mu^*\big),
$$
in the class of nonnegative definite matrices with $\sum_{i,j}\Gamma_{ij}=0$.

\end{rem}

In the two dimensional case $d=2$, the exact expression is known for $\M_\pi$ and one gets a precise formula for the
filter stability index
\begin{cor}
For  $d=2$ and  $\Delta h := h_1- h_2\ne 0$
\begin{equation}
\label{2d}
\gamma =
-(\lambda_{12}+\lambda_{21}) + \frac{(\Delta h)^2}{\sigma^2} \Big(-\frac{1}{2} +  \frac{\int_0^1 q(x)x(1-x)dx}{\int_0^1 q(x)dx}\Big)
\end{equation}
where
\begin{multline}
\label{q(x)}
q(x)=\frac{1}
{x^2(1-x)^2}\times\\
\exp\left(
-\frac{2\sigma^2\lambda_{21}}{(\Delta h)^2x(1-x)}
+\frac{2\sigma^2(\lambda_{12}-\lambda_{21})}{(\Delta h)^2}
\left(\log\frac{x}{1-x}+\frac{1}{1-x}\right)
\right).
\end{multline}
In particular,
\begin{equation}
\label{2d-mr1}
\gamma(\sigma)
 =
 -\frac{1}{2}\sigma^{-2}(\Delta h)^2+
 \big(1 +o (1)\big)\log\sigma^{-2}
\frac{2\lambda_{12}\lambda_{21}}{\lambda_{12}+\lambda_{21}},\quad \sigma\searrow 0
\end{equation}
and
\begin{multline}\label{2d-mr2}
\gamma(\sigma)
=  -(\lambda_{12}+\lambda_{21}) +
\sigma^{-2}\Big(-\frac{1}{2} \big(h^2_1+h^2_2\big)+
\big(\mu_2h_1+\mu_1 h_2\big)\big(\mu_1h_1+\mu_2 h_2\big)
\Big) +\\
 \sigma^{-4}\frac{(h_1-h_2)^2\mu^2_1\mu^2_2}{2(\lambda_{12}+\lambda_{21})} + o(\sigma^{-4}), \quad \sigma \nearrow \infty.
\end{multline}
\end{cor}
\begin{proof}
For $d=2$ the process $\rho_t\wedge \bar{\rho}_t$ is one dimensional (cf. \eqref{zeq} below) and
$Z_t:=\rho_t(1)\bar{\rho}_t(2)-\rho_t(2)\bar{\rho}_t(1)$ satisfies
\begin{align*}
dZ_t = -(\lambda_{12}+\lambda_{21})Z_t dt + \sigma^{-2}\big(h_1+h_2\big) Z_t dY_t+
\sigma^{-2}h_1h_2Z_t dt.
\end{align*}
If e.g. $Z_0>0$, then $Z_t>0$ for all $t\ge 0$ and by the It\^o formula
\begin{align*}
\lyap_1+\lyap_2=&\frac{1}{t}\log |Z_t|  = \frac{1}{t}\int_0^t \frac{1}{|Z_s|}dZ_s -\frac{1}{2}  \sigma^{-2}\big(h_1+h_2\big)^2 = \\
&-(\lambda_{12}+\lambda_{21}) + \sigma^{-2}\big(h_1+h_2\big)\Big( \frac{1}{t}\int_0^t h(X_s)ds +\frac{\sigma B_t}{t}\Big)+\\
&\sigma^{-2}h_1h_2
-\frac{1}{2}  \sigma^{-2}\big(h_1+h_2\big)^2\xrightarrow[]{t\to\infty} \\
& -(\lambda_{12}+\lambda_{21}) + \sigma^{-2}\Big(
\big(h_1+h_2\big)\mu(h)
-\frac{1}{2}  h^2_1-\frac{1}{2}h^2_2\Big).
\end{align*}
Note that in the two dimensional case the upper bound in \eqref{lyaplyap} is always attained.
The equation \eqref{W} is also one dimensional  and  $\pi_t:=\pi_t(1)$ satisfies
$$
d\pi_t = \big(\lambda_{21}-(\lambda_{12}+\lambda_{21})\pi_t\big)dt +\sigma^{-1}\pi_t(1-\pi_t)\big(h_1-h_2\big)d\bar{B}_t,
$$
where $\bar{B}$ is the innovation Brownian motion. The stationary probability distribution of $\pi_t$
has a density $q(x)$, solving the Kolmog\-orov-Fokker-Plank equation
$$
0=-\frac{\partial}{\partial x}\Big[\big(\lambda_{21}-(\lambda_{12}+\lambda_{21})x\big)q(x)\Big]+
\frac{\big(\Delta h\big)^2}{2\sigma^2} \frac{\partial^2}{\partial x^2}\Big[
x^2(1-x)^2q(x)
\Big],
$$
which has an explicit solution given by  \eqref{q(x)}.
The formula \eqref{2d} is nothing but \eqref{Lyap}, combined with the above expression for $\lyap_1+\lyap_2$ and the formula \eqref{L1},
where the explicit integration versus $q(x)$ appears.
The asymptotic \eqref{2d-mr1} is obtained by means of the expansion \eqref{refined}. The asymptotic \eqref{cr} and the expression
for $\lyap_1+\lyap_2$ give the first two terms in \eqref{2d-mr2}. The last term is obtained  via \eqref{crref}.\qed
\end{proof}

In fact both low and high signal-to-noise bounds \eqref{lyap} and \eqref{AZu} can be obtained by means of the same approach.
\begin{cor}
Assume that $X$ is ergodic, then for any $(\nu, \bar{\nu})$ the following global version of \eqref{lyap} holds:
\begin{equation}
\label{lowsnr}
\varlimsup_{\sigma\to \infty}\gamma(\sigma) \le  \gamma_{\max}(\Lambda).
\end{equation}
\end{cor}

\begin{proof}

In view of \eqref{Lyap} and \eqref{cr}, the claim \eqref{lowsnr} holds if
\begin{equation}
\label{toprove}
\varlimsup_{\sigma\to\infty}\lim_{t\to\infty}\frac{1}{t}\log |\rho^\sigma_t \wedge \bar{\rho}^\sigma_t|\le \gamma_{\max}(\Lambda),
\end{equation}
The process $Z^\sigma_t = \rho^\sigma_t \wedge \bar{\rho}^\sigma_t$ satisfies the linear equation
\begin{equation}
\label{zeq}
dZ^\sigma_t = \big(\Lambda^* Z^\sigma_t + Z^\sigma_t \Lambda\big)dt + \sigma^{-2} \big(HZ^\sigma_t +Z^\sigma_t H\big)dY_t +
\sigma^{-2}HZ^\sigma_t Hdt,
\end{equation}
subject to $Z^\sigma_0=\nu\wedge \bar{\nu}$, where $H:=\diag(h)$.

Let $U^\sigma(s,t)$ be the fundamental solution of \eqref{zeq}, i.e. the linear (random) operator such that
$
Z^\sigma_t =U^\sigma(s,t)\circ Z^\sigma_s
$, $\forall t\ge s\ge 0$.
Denote by $\Q^\nu$ and $\Q^\mu$  the probability measures, induced by $(X,Y)$ when $X_0\sim \nu$ and $X_0\sim \mu$
respectively. Since $\mu$ has strictly positive atoms, $\nu\ll\mu$ and  by the Markov property of the pair
$(X,Y)$ we have $\Q^\nu\ll \Q^\mu$. Then any event, which occurs $\Q^\mu$-a.s. occurs $\Q^\nu$-a.s. as well. In particular it is sufficient
to prove \eqref{toprove} for the stationary $\tilde X$ (note that $\rho_t$ and $\bar{\rho}_t$
are still the solutions of \eqref{zak} subject to $\nu$ and $\bar{\nu}$).

The idea of the proof is to ``sample'' the convergence in \eqref{toprove} on the subsequence $0,\tau,2\tau,...$ with a fixed $\tau>0$, to apply the LLN and
to study the obtained limit via $\lim_{\tau\to\infty}\lim_{\sigma\to \infty}$.
Define $\Pi^\sigma_t=Z^\sigma_t/|Z^\sigma_t|$
and let $\Dom$ be the set of all antisymmetric matrices equal to $\eta\wedge \eta'$ for some  $\eta, \eta'\in\simplex$ up to multiplication by a
nonzero constant. Clearly $\Dom$ is a closed set and the solution of \eqref{zeq} evolves in $\Dom$.
MET guarantees that  $\lim_{t\to\infty}\frac{1}{t}\log |Z^\sigma_t|$ exists $\P$-a.s. (see \cite{AZ97a}) and thus
for an arbitrary constant $\tau>0$
\begin{equation}\label{upp}
\begin{aligned}
\lim_{t\to\infty}\frac{1}{t}\log |Z^\sigma_t| =&
\lim_{k\to\infty}\frac{1}{k}\frac{1}{\tau}\log |Z^\sigma_{\tau k}|= \\
& \lim_{k\to\infty}\frac{1}{k}\sum_{m=1}^k\frac{1}{\tau}\log \big|U^\sigma(\tau(m-1),\tau m) \circ \Pi^\sigma_{\tau (m-1)}\big|\le \\
& \lim_{k\to\infty}\frac{1}{k}\sum_{m=1}^k\frac{1}{\tau}\log \max_{|v|= 1, v\in \Dom}\big|U^\sigma(\tau(m-1),\tau m) \circ v\big| =\\
& \frac{1}{\tau}\E \log \max_{|v|= 1, v\in \Dom}\big|U^\sigma(0, \tau) \circ v\big|
\end{aligned}
\end{equation}
The latter equality is due to the LLN, which holds since the summands form a stationary uniformly
integrable ergodic sequence. Indeed
$$u_m :=\log \max_{|v|= 1, v\in \Dom}\big|U^\sigma(\tau(m-1),\tau m) \circ v\big|$$
are measurable with respect to
$\F^X_{(m-1)\tau,m\tau}\vee \F^B_{(m-1)\tau,m\tau}$ and so stationarity
and ergodicity  are inherited from $X$, the increments of $B$ and their independence. Integrability follows from the Gaussian properties of $B$ and is verified
similarly to Theorem 1.5 in \cite{AZ97a}.

Observe that the solution of \eqref{zeq} converges uniformly on $[0,\tau]$ to the solution of
\begin{equation}
\label{Qeq}
\dot{Q}_t = \Lambda^* Q_t + Q_t\Lambda,\quad Q_0=\nu\vee \bar \nu
\end{equation}
as $\sigma\to\infty$, i.e.
$
\lim_{\sigma\to\infty} \sup_{t\in [0,\tau]}\big|Z^\sigma_t -Q_t \big|=0.
$
Consequently
\begin{equation}\label{opconv}
U^\sigma(0,\tau) \circ v\xrightarrow{\sigma\to \infty} V(0,\tau)\circ v,\quad \forall v\in \Dom
\end{equation}
where $V(s,t)$ is the fundamental solution of \eqref{Qeq}.
Notice that $V(0,t)\circ v=p_t\wedge q_t$, where
$p_t$ and $q_t$ solve $\dot{x}_t=\Lambda^* x_t$ subject to $p,q\in \simplex$ ($v=p\wedge q$). If $X$ is ergodic,
the zero eigenvalue of its transition rates matrix is simple and $\Lambda$ is a stability matrix on
$\{x\in\Real^d:\sum_{i=1}^d x_i =0\}$ (see e.g. \cite{Nor}):
$$
|p_t-q_t|= \big|e^{\Lambda^* t}(p-q)\big|\le c \exp\big(\gamma_{\max}(\Lambda)t\big),
$$
for a constant $c>0$.
Since
$
\frac{1}{2}|p_t\wedge q_t|\le |p_t-q_t|\le |p_t\wedge q_t|,
$
this implies
\begin{equation}\label{cl1}
\varlimsup_{\tau\to\infty}\frac{1}{\tau}\log \big|V(0,\tau)\circ v\big|\le \gamma_{\max}(\Lambda), \quad \forall v\in \Dom.
\end{equation}

Passing to the limits $\lim_{\tau\to\infty}\lim_{\sigma\to\infty}$ in \eqref{upp} and taking into account \eqref{opconv} and \eqref{cl1}, one gets
the required
$$
\varlimsup_{\sigma\to \infty}\lim_{t\to\infty}\frac{1}{t}\log |Z^\sigma_t|  \le \gamma_{\max}(\Lambda). \qed
$$

\end{proof}

\begin{cor}
Assume that $X$ is ergodic, then \eqref{AZu} holds.
\end{cor}
\begin{proof}
In view of Remark \ref{rem-1}, \eqref{AZu} holds if
\begin{equation}
\label{sum}
\lim_{\sigma\to 0}\sigma^2\big(\lambda_1(\sigma)+\lambda_2(\sigma)\big)\le \mu(h^2) - \frac{1}{2}\sum_{i=1}^d \mu_i\min_{j\ne i}\big(h(a_i)-h(a_j)\big)^2.
\end{equation}
Let $\zeta^\sigma_t := \rho_{t\sigma^2}$ and $\bar{\zeta}^\sigma_t := \bar{\rho}_{t\sigma^2}$ be the time scaled solutions of \eqref{zak}
subject to $\nu$ and $\bar{\nu}$ respectively. Then
$$
\lim_{t\to\infty} \frac{1}{t} \log |\rho^\sigma_t\wedge \bar{\rho}^\sigma_t| =
\frac{1}{\sigma^2} \lim_{t\to\infty} \frac{1}{t} \log |\zeta^\sigma_t \wedge \bar{\zeta}^\sigma_t|,
$$
where the limits exist by the Oseledec MET as mentioned before. The process $\zeta^\sigma_t$ (and $\bar{\zeta}^\sigma_t$)
satisfies the equation
$$
d\zeta^\sigma_t = \sigma^2 \Lambda^* \zeta^\sigma_t dt  + H \zeta^\sigma_t dW_t,
$$
where $W_t = \int_0^t h(X_{s\sigma^2})ds + \tilde B_t$ and $\tilde B_t := \sigma^{-1}B_{t\sigma^2}$ is a standard Brownian motion.
Consequently the process $R^\sigma_t = \zeta^\sigma_t \wedge \bar{\zeta}^\sigma_t$
satisfies the linear equation
$$
dR^\sigma_t = \sigma^2 \big(\Lambda^* R^\sigma_t + R^\sigma_t \Lambda\big)dt +  \big(HR^\sigma_t +R^\sigma_t H\big)dW_t +
H R^\sigma_t Hdt,
$$
subject to $R^\sigma_0=\nu\wedge \bar{\nu}$. In the componentwise notation the latter reads
\begin{multline}
\label{wedge-sde}
dR^\sigma_{km}(t) = \sigma^2\Big(\sum_{j\ne k}\lambda_{jm}R^\sigma_{kj}(t)+\sum_{j\ne m}\lambda_{jk}R^\sigma_{jm}(t)\Big)dt + \\
(h_m+h_k)R^\sigma_{km}(t) dW_t
+h_mh_k R^\sigma_{km}(t) dt, \quad k\ne m.
\end{multline}
Let $U^\sigma(s,t)$ denote the fundamental solution of \eqref{wedge-sde} (cf. \eqref{zeq}), i.e. a tensor whose entries $U^\sigma_{ij,km}(s,t)$ are
$\F^W_{s,t}$-measurable and
\begin{equation}\label{Qsig}
R^\sigma_{km}(t) = \Big[U^\sigma(s,t)\circ R^\sigma_s\Big]_{km}=\sum_{i\ne j} U^\sigma_{ij,km}(s,t)R^\sigma_{ij}(s), \quad k\ne m.
\end{equation}
For fixed $i\ne j$, $U^\sigma_{ij,km}(s,t)$  are the entries of the  matrix generated by \eqref{wedge-sde},
subject to $$U^\sigma_{ij,km}(s,s)=\delta_{ij}^{km}:=\begin{cases}
1, & i=k, j=m\\
0,& \text{otherwise}.
\end{cases}$$
Define
\begin{multline*}
\psi^\sigma_{km}(s,t) := \exp\big(\sigma^2(\lambda_{kk}+\lambda_{mm})(t-s)+\\
(h_k+h_m)(W_t-W_s)+h_kh_m(t-s)-\frac{1}{2}(h_k+h_m)^2(t-s)\big).
\end{multline*}
Then for $i\ne j$ and $k\ne m$
\begin{multline}\label{repr}
U^\sigma_{ij,km}(s,t) = \psi^\sigma_{km}(s,t)\delta_{ij}^{km} + \\ \sigma^2 \psi^\sigma_{km}(s,t)\int_s^t \big(\psi^\sigma_{km}(s,u)\big)^{-1}
\sum_{q\ne k,m}\Big(\lambda_{qm}U^\sigma_{ij,kq}(s,u)+\lambda_{qk}U^\sigma_{ij, qm}(s,u)\Big)
du.
\end{multline}
Since the off-diagonal entries of $\Lambda$ are nonnegative, the latter implies that $U^\sigma(s,t)$ is a positive operator.
Let $V^\sigma_{ij,km}(s,t)=U^\sigma_{ij,km}(s,t)/\sum_{k\ne m}U^\sigma_{ij,km}(s,t)$,  then
\begin{align*}
&\log \Big(\sum_{k\ne m}U^\sigma_{ij,km}(s,t)\Big) =
 \int_s^t
\sigma^2\sum_{k\ne m}\Big(\sum_{\ell\ne k}\lambda_{\ell m}
V^\sigma_{ij,km}(s,r)
+\sum_{\ell\ne m}\lambda_{\ell k}V^\sigma_{ij,\ell m}(s,r)\Big)dr  \\
& \hskip 0.5in +\int_s^t \sum_{k\ne m}(h_m+h_k)V^\sigma_{ij,km}(s,r) dW_r
+\int_s^t \sum_{k\ne m} h_mh_k V^\sigma_{ij,km}(s,r) dr -\\
& \hskip 0.5in \frac{1}{2} \int_s^t \Big(\sum_{k\ne m}(h_m+h_k)V^\sigma_{ij,km}(s,r)\Big)^2dr\le\\
& \hskip 0.5in \big( 2\sigma^2 d\lambda_{\max}+ 5h^2_{\max}\big)(t-s) + \int_s^t \sum_{k\ne m}(h_m+h_k)V^\sigma_{ij,km}(s,r)d\tilde B_r
\end{align*}
where $\lambda_{\max}=\max_{i\ne j}\lambda_{ij}$ and $h_{\max}=\max_i |h_i|$.  The latter and \eqref{repr} gives the
following estimate
$$
U^\sigma_{ij,km}(s,t) \le  \psi^\sigma_{km}(s,t)\delta_{ij}^{km} + \sigma^2 C_1 e^{C_2(t-s)} \int_s^t \exp\left(\int_s^u \beta^\sigma_{ij,km}(r) d\tilde B_r\right)du \\
$$
with some constants $C_1$ and $C_2$ and bounded processes $\beta^\sigma_{ij,km}(t)$. Moreover
$$
\varphi^\sigma_{ij,km}(s,t) := C_1 e^{C_2(t-s)} \int_s^t \exp\left(\int_s^u \beta^\sigma_{ij,km}(r) d\tilde B_r\right)du
$$
are uniformly integrable in $\sigma$, since  bounded $\beta^\sigma_{ij,km}(r)$ trivially satisfy the Novikov condition
$
\E\exp\left(\frac{1}{2} \int_s^t \big(\beta^\sigma_{ij,km}(r)\big)^2 dr\right)<\infty
$ (see Theorem 6.1 in \cite{LS}).
Define $\Pi^\sigma(t)=R^\sigma(t)/|R^\sigma(t)|$ and fix a $\tau>0$, then for $t\ge \tau$ (see \eqref{Qsig})
\begin{align*}
& \log|R^\sigma(t)| = \log|R^\sigma(t-\tau)| +
\log \big|\sum_{\substack{i\ne j\\ k\ne m}} U^\sigma_{ij,km}(t-\tau,t)\Pi^\sigma_{ij}(t-\tau)\big| \le \\
&
\log|R^\sigma(t-\tau)| +
\log \max_{i\ne j}\sum_{k\ne m} U^\sigma_{ij,km}(t-\tau,t)\le \\
&
\log|R^\sigma(t-\tau)| +
\log \Big(\max_{i\ne j} \psi^\sigma_{ij}(t-\tau,t)
+\sigma^2
\max_{i\ne j}\sum_{k\ne m}\varphi^\sigma_{ij,km}(t-\tau,t)
\Big).
\end{align*}
As was mentioned before, it is enough to establish \eqref{sum} for the stationary chain $\tilde X$.
Recall that by  MET the limit $\lim_{t\to\infty}\frac{1}{t}\log |R^\sigma_t|$ exists $\P$-a.s. Then for any $\nu, \bar{\nu}\in \simplex$ and $\tau>0$
\begin{align*}
& \lim_{t\to\infty}\frac{1}{t}\log |R^\sigma_t| = \lim_{n\to\infty}\frac{1}{n\tau}\log |R^\sigma_{n\tau}|
\le \\
& \varlimsup_{n\to\infty}\frac{1}{n} \sum_{\ell=1}^n
\frac{1}{\tau}\log
\Big(\max_{i\ne j} \psi^\sigma_{ij}\big((\ell-1)\tau,\ell\tau\big) +
\sigma^2
\max_{i\ne j}\sum_{k\ne m}\varphi^\sigma_{ij,km}\big((\ell-1)\tau,\ell \tau\big)\Big)=\\
& \frac{1}{\tau}\E
\log
\Big(\max_{i\ne j} \psi^\sigma_{ij}\big(0,\tau\big) +
\sigma^2
\max_{i\ne j}\sum_{k\ne m}\varphi^\sigma_{ij,km}(0,\tau)
\Big).
\end{align*}
The latter equality holds by LLN, since being measurable functionals of the increments of $\tilde X$ and $\tilde B$ on the intervals $[(\ell-1)\tau,\ell\tau)$,
both $\psi^\sigma_{ij}\big((\ell-1)\tau,\ell\tau\big)$ and $\varphi^\sigma_{ij,km}\big((\ell-1)\tau,\ell \tau\big)$, $\ell\ge 1$ form stationary sequences with
finite expectations.
On the set $A^\sigma_\tau:=\{\tilde X_0=\tilde X_{u\sigma^2}, u\le \tau\}$ we have
\begin{align*}
&\psi^\sigma_{km}(0,\tau)=\exp\big(\sigma^2(\lambda_{kk}+\lambda_{mm})\tau+(h_k+h_m)W_\tau+\\
 & \hskip 3.0in h_kh_m\tau-\frac{1}{2}(h_k+h_m)^2\tau\big) =\\
&\exp\big(h_kh_m\tau-\frac{1}{2}(h_k+h_m)^2\tau+(h_k+h_m)h(\tilde X_0)\tau +\\
& \hskip 2.48in(h_k+h_m)\tilde B_\tau+\sigma^2(\lambda_{kk}+\lambda_{mm})\tau\big) =\\
& \exp\big(h^2(\tilde X_0)\tau-\frac{1}{2}(h_k-h(\tilde X_0))^2\tau-\frac{1}{2}(h_m-h(\tilde X_0))^2\tau +\\
& \hskip 2.48in (h_k+h_m)\tilde B_\tau+\sigma^2(\lambda_{kk}+\lambda_{mm})\tau\big)\le \\
& \exp\big(h^2(\tilde X_0)\tau-\frac{1}{2}\min_{a_m\ne \tilde X_0}(h_m-h(\tilde X_0))^2\tau+
2h_{\max}|\tilde B_\tau|+2d\sigma^2 \lambda_{\max} \tau\big).
\end{align*}
Since the process $\tilde X_{t\sigma^2}$ is a slow chain with the generator $\sigma^2\Lambda$, $\lim_{\sigma\to 0}\P_s(A^\sigma_\tau)=1$
for any fixed $\tau>0$.
This gives the following estimate
\begin{align*}
& \lim_{t\to\infty}\frac{1}{t}\log |R^\sigma_t| \le \\
& \hskip 0.1in \frac{1}{\tau}\E\one{A^\sigma_\tau}\big(h^2(\tilde X_0)\tau-\frac{1}{2}\min_{a_m\ne \tilde X_0}(h_m-h(\tilde X_0))^2\tau+
2h_{\max}|\tilde B_\tau|+2d\sigma^2 \lambda_{\max} \tau\big) + \\
& \hskip 0.1in \frac{1}{\tau}\E\one{\Omega\backslash A^\sigma_\tau}
\log
\big(\max_{i\ne j} \psi^\sigma_{ij}(0,\tau) +
\sigma^2\max_{i\ne j}\sum_{k\ne m}\varphi^\sigma_{ij,km}(0,\tau)\big)
\xrightarrow{\sigma\to 0}\\
&
\hskip 0.1in
\E \Big(h^2(\tilde X_0)-\frac{1}{2}\min_{a_m\ne \tilde X_0}(h_m-h(\tilde X_0))^2
+\frac{2h_{\max}|\tilde B_\tau|}{\tau}\Big)
\end{align*}
which implies \eqref{sum} as $\tau\to \infty$.\qed
\end{proof}

\section{Conclusions}

The Wonham and Kalman-Bucy filters are particular instances of the general filtering equation, for which
finite dimensional realizations are known. Consequently both are of considerable practical interest in various applications (see e.g. \cite{EAM}).
While the stability and ergodic properties for the linear Kalman-Bucy filter has been studied and understood since 60's, the analogous
theory for the Wonham filter is less developed and in fact was addressed only a decade ago. In this paper we established
certain ergodic properties of the signal/filtering pair, which are crucial for applying the classic Lyapunov exponents technique for
SDEs (\cite{Kh}) in the filtering context. The latter allows to derive refined formulae for the Lyapunov exponents of the Zakai equation and
simplify derivation of certain known bounds on the filter stability index. In particular case of the binary signal, a complete characterization
of the filter stability is obtained.

\vskip 0.1in
\noindent
{\bf Acknowledgements}. I am grateful to Zvi Artstein and Ramon van Handel for their suggestions regarding this paper.
The paper presentation has been significantly improved by the comments of anonymous referee.


\begin{thebibliography}{}
\bibitem{Arnold} L. Arnold, Random dynamical systems. Springer Monographs in Mathematics. Springer-Verlag, Berlin, 1998

\bibitem{AZ97a} R. Atar, O. Zeitouni,  Lyapunov exponents for finite state nonlinear filtering. SIAM J. Control Optim. 35 (1997), no. 1, 36--55.

\bibitem{AZ97b} R. Atar, O. Zeitouni, Exponential stability for nonlinear filtering, Ann. Inst. H. Poincare Probab. Statist. 33 (1997), no. 6, 697--725

\bibitem{BO} A.Budhiraja, D.Ocone, Exponential stability in discrete-time filtering for non-ergodic signals. Stochastic Process. Appl. 82 (1999), no. 2, 245--257

\bibitem{Ch3} P. Chigansky, On filtering of Markov chains in strong noise, to appear in IEEE Trans. Inf. Theory, preprint \verb"math.PR/0508446" at \verb"www.arxiv.org"

\bibitem{DZ} B.~Delyon, O.~Zeitouni,  Lyapunov exponents for filtering problem, in Applied Stochastic Analysis,  Davis, M. H. A. and
Elliot R. J. eds., Gordon \& Breach, New York, 1991, pp.~511-521.

\bibitem{EAM} R.J.Elliott, L.Aggoun, J.B.Moore, Hidden Markov models. Estimation and control. Applications of Mathematics
(New York), 29. Springer-Verlag, New York, 1995

\bibitem{BK} A. Budhiraja, H.J.Kushner,  Approximation and limit results for nonlinear filters over an infinite time interval. SIAM J. Control Optim. 37 (1999), no. 6, 1946--1979

\bibitem{BCL} P. Baxendale, P. Chigansky, R.Liptser,  Asymptotic stability of the Wonham filter: ergodic and nonergodic signals,
SIAM J. Control Optim. 43 (2004), no. 2, 643--669.

\bibitem{DG} P. Del Moral, A. Guionnet, On the stability of interacting processes with applications to filtering and genetic algorithms.
Ann. Inst. H. Poincare Probab. Statist. 37 (2001), no. 2, 155--194

\bibitem{Gol} G. Golubev, On filtering for a Hidden Markov Chain under square perfromance critetion,
Problems of Information Transmissionm Vol. 36, No. 3, 2000, pp. 213-219

\bibitem{Kh} R. Khasminskii, Stochastic stability of differential equations, Imprint Alphen aan den Rijn, The Netherlands, Sijthoff \& Noordhoff, 1980
(this note cites the Russian original, Moscow, 1969)

\bibitem{KZ} R. Khasminskii, O. Zeitouni, Asymptotic filtering for finite state Markov chains, Stoch.
Processes and Appl., 1996, vol. 63, pp. 1-10.

\bibitem{Kuflows} H. Kunita, Stochastic differential equations and stochastic flows of diffeomorphisms. Ecole d'ete de probabilites de Saint-Flour,
XII--1982, 143--303, Lecture Notes in Math., 1097, Springer, Berlin, 1984

\bibitem{Ku} H. Kunita, Asymptotic behavior of the nonlinear filtering errors of Markov processes. J. Multivariate Anal. 1 (1971), 365--393

\bibitem{LO} F.LeGland, N. Oudjane, A robustification approach to stability and to uniform particle approximation of nonlinear filters:
the example of pseudo-mixing signals. Stochastic Process. Appl. 106 (2003), no. 2, 279--316

\bibitem{LS} R.Lipster, A.Shiryaev, Statistics of random processes: theory and applications, I,II, Springer-Verlag, 2nd Ed., 2001

\bibitem{Nor} J.R.Norris, Markov chains, Cambridge Series in Statistical and Probabilistic Mathematics, 2. Cambridge University Press,
Cambridge, 1998

\bibitem{RS} B.Rozovskii, A. Shiryaev, On infinite systems of stochastic differential equations that arise in the theory of optimal
nonlinear filtering, Teor. Verojatnost. i Primenen. 17 (1972), 228--237


\bibitem{Var} S.R.S. Varadhan, Probability theory, Courant Lecture Notes in Mathematics, 7. New York University, Courant Institute of
 Mathematical Sciences, New York; American Mathematical Society, Providence, RI, 2001

\bibitem{W} W. M. Wonham, Some applications of stochastic differential equations to optimal nonlinear filtering.
J. Soc. Indust. Appl. Math. Ser. A Control 2 347--369 (1965).

\end{thebibliography}
\end{document}